\documentclass[12pt,twoside,a4paper]{article}
\usepackage{amsmath,amssymb}
\usepackage{graphicx}
\usepackage{fullpage}
\usepackage{pslatex}
\usepackage{verbatim}
\usepackage{url}

\newcommand{\Odip}[2]{\mathcal{O}_{#1}\!\left(#2\right)\mathchoice{\!}{}{}{}}
\newcommand{\Odi}[1]{\Odip{}{#1}}

\newcommand{\odip}[2]{{o}_{#1}\!\left(#2\right)\mathchoice{\!}{}{}{}}
\newcommand{\odi}[1]{\odip{}{#1}}

\newcommand{\dx}{\mathrm{d}}

\newcommand{\A}{\mathcal{A}}
\newcommand{\B}{\mathcal{B}}
\newcommand{\C}{\mathcal{C}}

\renewcommand{\P}{\mathcal{P}}
\newcommand{\PP}{\mathfrak{P}}

\newcommand{\RR}{\mathcal{R}}

\newcommand{\R}{\mathbb{R}}
\newcommand{\Z}{\mathbb{Z}}

\newcommand{\singseries}{\mathfrak{S}}

\newcommand{\Odiprimep}[2]{\mathcal{O'}_{#1}\!\Bigl(#2\Bigr)
				\mathchoice{\!}{}{}{}}

\newcommand{\Stirling}[2]{\genfrac{\{}{\}}{0pt}{}{#1}{#2}}

\newtheorem{Theorem}{Theorem}
\newtheorem{Corollary}{Corollary}
\newtheorem{Lemma}{Lemma}

\newenvironment{Proof}[1][Proof]{\par\noindent\textbf{#1.}~}
  {\hfill$\square$\smallskip\par}

\begin{document}

\title{Prime numbers in logarithmic intervals}
\author{D.~BAZZANELLA, A.~LANGUASCO, A.~ZACCAGNINI}
\date{September 17, 2008}

\maketitle

\begin{abstract}
Let $X$ be a large parameter.
We will first give a new estimate for the integral moments
of primes in short intervals of the type $(p,p+h]$,
where $p\leq X$ is a prime number and $h=\odi{X}$. Then 
we will apply this to 
prove that for every $\lambda>1/2$ there exists a positive proportion
of primes $p\leq X$ such that the interval $(p,p+ \lambda\log X]$ contains 
at least a prime number. As a consequence we improve Cheer and Goldston's
result on the size of real numbers $\lambda>1$ with the property that
there is a positive proportion of integers $m\leq X$ such that
the interval $(m,m+ \lambda\log X]$ 
contains no primes. We also prove other results
concerning the moments of the gaps between consecutive primes
and about the positive  proportion of integers $m\leq X$ such that
the interval $(m,m+ \lambda\log X]$ 
contains at least a prime number. 
The last application of these techniques are two theorems
(the first one unconditional and the second one in which we assume the 
validity of the Riemann Hypothesis and of a form of the Montgomery
pair correlation conjecture) on
the positive proportion of primes $p\leq X$ such that the interval 
$(p,p+ \lambda\log X]$ contains no primes.
\end{abstract}

\section{Introduction}

Let $X$ be a large parameter, $\PP$ be the set of primes and $\lambda$
be a positive real number.
This paper is devoted to study the distribution of primes in short
intervals: in particular we will give lower bounds
for the proportion of positive integers $m\leq X$, or $p\in \PP$
and $p\leq X$,
such that the intervals 
$(m,m+ \lambda\log X]$ or $(p,p+ \lambda\log X]$
contain or do not contain a prime number.

Many mathematicians studied the distribution of primes in short intervals;  
here we just recall some fundamental papers on this topic. Several  results are formulated using the quantity
\[
  E
  =
  \liminf_{i\to+\infty} \frac{p_{i+1}-p_i}{\log p_i}.
\]
In 1926 Hardy and Littlewood \cite{HardyLittlewood1926}, assuming the validity 
of the Generalized Riemann Hypothesis, gave the first non-trivial estimate
$E\leq 2/3$. 
In 1940, Erd{\H o}s \cite{Erdos1940} proved unconditionally that $E<1$ and in 
1966 Bombieri and Davenport \cite{BombieriD1966} improved this result to
$E\leq 0.46650\dotsc$.
In 1972-73, Huxley \cite{Huxley1972,Huxley1973}, using
a new set of weights, was able to reach $E\leq 0.44254\dotsc$.
For all these results, a suitable modification of the argument
can lead to a positive proportion result on the cardinality of the
integers $m\leq X$ such that $(m,m+\lambda \log X]$ contains at
least a prime number, for any fixed $\lambda$ larger than the given
bound for $E$. In 1986 Maier used his 
matrix-method \cite{Maier1988}
to obtain $E\leq 0.2486\dotsc$ but this method gave 
no positive proportion results.

It was just in 2005 that Goldston, Pintz and
Y{\i}ld{\i}r{\i}m \cite{GoldstonPY2005} obtained 
that  $E=0$ solving a long-standing conjecture. 
In fact they proved quite a stronger result:
\[
\textrm{there are infinitely many}\ i\ \textrm{such that} \
p_{i+1}-p_i \ll \sqrt{\log p_i}\ (\log \log p_i)^2.
\]
Unfortunately, it seems that
this wonderful new technique gives no positive proportion
results.

In 1987 Cheer and Goldston \cite{CheerG1987b,CheerG1987a} 
used the integral moments of primes over integers and a refinement
of Erd{\H o}s's technique to prove results of this kind.
We also recall that Goldston and Y{\i}ld{\i}r{\i}m \cite{GoldstonY2007}, 
in a paper published in 2007 but that was developed before 
\cite{GoldstonPY2005} appeared, were able to obtain a new proof 
of the inequality $E\leq 1/4$ with a method which gives also a 
positive proportion result.

Here we use a new result on integral moments of primes over primes
(see Theorem \ref{Thm-huge-costant}), to prove that there exists a
positive proportion of primes $p\leq X$ such that the interval
$(p,p+ \lambda\log X]$ contains at least a prime number
for every $\lambda>1/2$, see Theorem \ref{Thm-A-Cardinality}. 
Even if the uniformity in $\lambda$ is weaker than the one proved in 
Theorem 1 of Goldston and Y{\i}ld{\i}r{\i}m \cite{GoldstonY2007}
($\lambda>1/4$), here we obtain an evaluation of implicit constant
which will be useful in the consequences.
The first of them is Theorem \ref{Thm-B-Cardinality} in which we
improve Cheer and Goldston's \cite{CheerG1987b} result on the size of
$\lambda>1$ such that there is a positive
proportion of integers $m\leq X$ such that
the interval $(m,m+ \lambda\log X]$ contains no primes.
The second consequence of Theorem \ref{Thm-A-Cardinality}
is a result concerning the moments of the gaps between 
consecutive primes
(see Theorem \ref{alpha-lower-bound}).

Theorem \ref{A1-lower-bound} is
about the positive proportion
of primes $p\leq X$ such that the interval $(p,p+ \lambda\log X]$,
where $\lambda$ is ``small'', contains 
no primes and its Corollary \ref{B-lower-bound} concerns
the positive proportion of integers $m\leq X$ such that
the interval $(m,m+ \lambda\log X]$ 
contains at least a prime number. 
Our last result (Theorem \ref{A1-cond-lower-bound}) slightly refines a conditional
theorem of  Cheer and Goldston \cite{CheerG1987b} on the 
positive proportion of primes $p\leq X$ such that the interval 
$(p,p+ \lambda\log X]$, where $\lambda>0$, contains 
no primes.

To be more precise in describing our results we need 
now to give some notation and definition.
Let $1\leq n \leq X/2$ be an integer.
We define the twin-prime counting functions as follows
\begin{equation} 
\label{Z-def}
Z(X; 2n)
=
\sum_{p\leq X} 
\sum_{\substack{p'\leq X \\ p' -p =2n}} \log p \log p'
\quad
\textrm{and}
\quad
Z_1(X; 2n)
=
\sum_{p\leq X} 
\sum_{\substack{p'\leq X \\ p' -p =2n}} 1
\end{equation}
where $p,p' \in \PP$.
Moreover we will write 
\begin{equation}
\label{singseries-def}
  \singseries(n)
  =
  2 c_0 \prod_{\substack{p \mid n \\ p > 2}} \frac{p - 1}{p - 2},
  \qquad\textrm{where}\qquad
  2 c_0
  =
  2 \prod_{p > 2} \Big( 1 - \frac{1}{(p-1)^2} \Big)
\end{equation}
to denote the singular series for this problem
and the twin-prime constant.
Letting  $1 < K \leq X$ be a real number, we define
\begin{equation} 
\label{A(K)-def}
\A(K)=
\{
p \leq X,
\ p \ \textrm{prime, such that
there exists a prime}\
p' \ \textrm{with} \ 
0< p'-p\leq K
\},
\end{equation}
\[
\A_1(K)
=  
\{
p \leq X,
\ p \ \textrm{prime} 
\}
\setminus \A(K),
\]
\[
\B(K)=
\{m \in [1,X]\cap \Z \
\textrm{such that
there exists a prime in}\
(m, m+K]
\}
\]
and
\begin{equation} 
\label{B1(K)-def}
\B_1(K) 
 =  
 [1,X] \setminus \B(K).
\end{equation}  
Let moreover 
\begin{equation}
\label{P-def}
\P_k(y)
=
\sum_{r=1}^{k} 
\Stirling{k}{r} 2^r r! 
y^r
\end{equation}
be a polynomial in $y$ where
$\Stirling{k}{r}$ denotes the Stirling number of second 
type defined as the number of ways to partition a set with $k$ elements 
into non-empty subsets having $r$ elements each 
(without counting the order of the subsets).

Recalling that $\pi(u)$ is the number of primes 
up to $u$ and that $\psi(u)= \sum_{n \leq u} \Lambda(n)$, where
$\Lambda(n)$ is the von Mangoldt function, we are now ready
to state the following result  about integral moments of primes
over primes in 
short intervals.
In what follows we will also denote by $\epsilon$ a small positive 
constant, not necessarily the same at each occurrence,
and by $\omega\geq1$ a parameter that will be useful in the applications.
\begin{Theorem}
\label{Thm-huge-costant}
Let $\epsilon>0$, $\omega>1$, and $h\in \R$, $f(X) \leq  h\leq X^{1-\epsilon}$,
where $f(X) \to +\infty$ arbitrarily slowly as $X\to +\infty$.
Let further $k \geq 2$ be an integer. 
Then
\[
\sum_{p \leq X}
(\psi(p+h) -\psi(p))^k
\leq
\Bigl(
\P_{k+1}
\Bigl(
\frac{\omega\, h}{\log X}
\Bigr)
+ \epsilon
\Bigr)
\frac{X}{h(\omega-1)}\log^{k} X,
\]
where $\P_k(y)$ is defined in \eqref{P-def}.
\end{Theorem}
The limitation to $\omega>1$ in Theorem \ref{Thm-huge-costant} 
and in the following applications arises from Lemma \ref{monotonicity} 
below.
We are mainly interested to the case $h= \lambda \log X$ where 
$\lambda>0$ is a constant. 
Letting
\begin{equation}
\label{R-def}
\RR_{\,\, k,\omega}(\lambda)
=
\frac{
\P_{k}(\omega \lambda)}{(\omega-1)\lambda}
=
\frac{1}{\omega-1}
\sum_{r=1}^{k} 
\Stirling{k}{r} r! 
2^r \omega^r \lambda^{r-1},
\end{equation}
we have the  
\begin{Corollary}
\label{Corollary-log}
Let $\epsilon>0$, $\omega>1$ and $k \geq 2$ be an integer. 
Let further $\lambda>0$ be a fixed constant. Then
\[
\sum_{p \leq X}
(\psi(p+\lambda \log X) -\psi(p))^k
\leq
\Bigl(
\RR_{\,\, k+1,\omega}(\lambda)+\epsilon
\Bigr)
X\log^{k-1} X
\]
and
\[
\sum_{p \leq X}
(\pi(p+\lambda \log X) -\pi(p))^k
\leq
\Bigl(
\RR_{\,\, k+1,\omega}(\lambda)+\epsilon
\Bigr)
\frac{X}{\log X}.
\]
\end{Corollary}
Denoting by $\vert \C \vert$ the cardinality of a given set $\C$,
we can now state our result on $\vert \A(K) \vert$
when $K$ is about $\log X$.
\begin{Theorem}
\label{Thm-A-Cardinality}
Let $\epsilon>0$, 
$X$ be a large parameter and $\A(K)$ be defined as in \eqref{A(K)-def}.
Let further $\lambda>1/2$ be a fixed constant.  We have that 
\[
\vert \A(\lambda \log X) \vert
\geq
( c_1(\lambda)-\epsilon) 
\frac{X}{\log  X},
\]
where
$c_1(\lambda)
= 
\sup_{\ell \in \Z;\  \ell \geq 2} 
\sup_{\omega>1} 
\Delta_{\ell,\omega}(\lambda) $ and
\begin{equation}
\label{Delta-def}
\Delta_{\ell,\omega}(\lambda)
=
\frac
{(\lambda/2 -1/4)^{\ell/(\ell-1)}}
{\RR_{\ \ell+1,\omega}(\lambda)^{1/(\ell-1)}}.
\end{equation}
\end{Theorem}
Theorem \ref{Thm-A-Cardinality} means that there is a positive proportion 
of primes $p\leq X$ such that the interval $(p,p+\lambda \log X]$, with $\lambda>1/2$, contains at least a prime number.  Our uniformity in $\lambda$ is weaker
than the one in Theorem 1 of Goldston and Y{\i}ld{\i}r{\i}m \cite{GoldstonY2007}  
($\lambda>1/4$) but there they gave no evaluation of the implicit constant.
Since in the following applications we will need this,
we have to use our weaker Theorem  \ref{Thm-A-Cardinality}.
 
Assuming a suitable form of the $k$-tuple conjecture,
see, \emph{e.g.}, equation \eqref{ktuple-conjecture} below, 
it is clear that equation \eqref{sq-brackets-estim} below holds 
for every positive $\lambda$ and with the factor $\lambda/2-1/4$ 
replaced by $\lambda/2$.
Hence in this case we can replace,
in the statement of Theorem \ref{Thm-A-Cardinality},
the condition $\lambda>1/2$ 
with $\lambda>0$ 
and
$\Delta_{\ell,\omega}(\lambda)$
with 
\begin{equation}
\label{Delta-tilde-def}
\widetilde{\Delta}_{\ell,\omega}(\lambda)
=
\frac
{(\lambda/2)^{\ell/(\ell-1)}}
{\widetilde{\RR}_{\ \ell+1,\omega}(\lambda)^{1/(\ell-1)}} ,
\end{equation}
where
\begin{equation}
\label{R-tilde-def}
\widetilde{\RR}_{\ \ell,\omega}(\lambda)
=
\frac{1}{\omega-1}
\sum_{r=1}^{\ell} 
\Stirling{\ell}{r}
\omega^r \lambda^{r-1}.
\end{equation}
In fact a simpler form of the constant in Theorem \ref{Thm-A-Cardinality} 
can be proved using $\omega=2$ and the Cauchy-Schwarz inequality 
instead of the H\"older inequality.
But numerical computations proved, at least for $\lambda\in (1/2,2]$, 
that the largest constants $\Delta_{\ell,\omega}(\lambda)$ 
and $\widetilde{\Delta}_{\ell,\omega}(\lambda)$ are  obtained 
with $\omega\approx 5503/5000$ and
$\ell=11$ in the unconditional case, and for $\omega\approx 5939/5000$
and  $\ell= 10$ in the conditional case;
moreover in the following application  a different optimization is needed
and the form in \eqref{Delta-def} is a more flexible one and leads to better  
final results. 
To write some numerical values, we have for $\lambda\in (1/2,2]$ that
the largest $\Delta_{\ell,\omega}(\lambda)$ is about $0.01266456\dotsc$ 
%
%
while the largest $\widetilde{\Delta}_{\ell,\omega}(\lambda)$ is about 
$0.11604228\dotsc$. 

As an application of Theorem~\ref{Thm-A-Cardinality} we have a result
about the set $\B_1(K)$.
We improve the estimates of Theorem~3 in Cheer and Goldston
\cite{CheerG1987b}. We also remark that, even if in \cite{CheerG1987b}
Cheer and Goldston used  $p_i\in (X,2X]$ while we are working
with $p_i\in (0,X]$, we still can compare the constants involved
since the estimates have a good dependence on $X$.
\begin{Theorem}
\label{Thm-B-Cardinality}
Let $\epsilon>0$, $X$ be a large parameter and $\B_1(K)$ be defined as
in \eqref{B1(K)-def}. Then there exists $\lambda>1$ such that
\begin{equation}
\label{B1-lowerbound}
  \vert \B_1(\lambda \log X) \vert
  \geq
  \bigl( c_2(B,\lambda) - \epsilon \bigr)
  X,
\end{equation}
where
\[
c_2(B,\lambda)
=
\sup_{\ell \in \Z;\  \ell \geq 2} \sup_{\omega>1}
\sup_{\nu\in(1/2,1-1/(2B))}
\
\frac
{B
\Bigl(
1-\lambda
+
(1-\Delta_{\ell,\omega}(\nu))^2/(2 B)
+
(\lambda-\nu)\Delta_{\ell,\omega}(\nu)
\Bigr)^2
}
{2(1-\Delta_{\ell,\omega}(\nu))^2}
\]
is a  positive constant,
$B$ is defined in Lemma \ref{BD-Thm2} below 
and $\Delta_{\ell,\omega}(\nu)$ is defined
in \eqref{Delta-def}. Moreover we also have
\begin{equation}
\label{dn-squared-lowerbound}
\sum_{p_i \leq X}
(p_{i+1}-p_i)^2
\geq
\Bigl(
1+ \frac{1}{12B^2}
+
c_3(B)
-\epsilon
\Bigr)
X\log X,
\end{equation}
where
\[
c_3(B)
=
\sup_{\ell \in \Z;\  \ell \geq 2} \sup_{\omega>1}
\sup_{\nu\in(1/2,1-1/(2B))}
\Bigr[
\frac{(1-\Delta_{\ell,\omega}(\nu))^3}{3B^2}
-
\frac{B}{3}
\Bigl(
\frac{\Delta_{\ell,\omega}(\nu)(\nu-1)}
{1-\Delta_{\ell,\omega}(\nu)}
+\frac{2-\Delta_{\ell,\omega}(\nu)}
{2B}
\Bigr)^3
\Bigr]
\]
is a  positive constant and 
$B$ 
and $\Delta_{\ell,\omega}(\nu)$ are as before. 
\end{Theorem}
The best $\lambda>1$ we are able to obtain in the previous statement is
\begin{equation}
\label{improved-lambda}
\begin{split}
\lambda
&
=
\sup_{\ell \in \Z;\  \ell \geq 2} \sup_{\omega>1}
\sup_{\nu\in(1/2,1-1/(2B))}
\Bigl(
1-\Delta_{\ell,\omega}(\nu)
\Bigr)^{-1}
\Bigl(
1
-
\nu\Delta_{\ell,\omega}(\nu)
+
\frac
{(1-\Delta_{\ell,\omega}(\nu))^2}
{2 B}
\Bigr)
-\epsilon
\\
&\ge
1.145358\dotsc
,
\end{split}
\end{equation}
%
where the numerical value is obtained using
the estimate of Fouvry and Grupp \cite{FouvryG1986} for $B=3.454$,
$\nu \approx 0.666856$, $\ell=12$ and $\omega\approx 2491/2250$ in 
\eqref{Delta-def}. 
Since it is not completely clear how to optimize the estimates in this
theorem, it is possible that the numerical values written here and later can be
further improved.
We remark that Cheer and Goldston  \cite{CheerG1987b}, in their Theorem~3,
proved that $\lambda = 1+1/(2B)$ is allowed in \eqref{B1-lowerbound}. 
Using the estimate for $B$ mentioned above,
this leads to $\lambda= 1.144759\dotsc$.

Moreover, we remark that the constant $c_2(B,\lambda)$ 
is larger than the one in eq.~(3.3) of \cite{CheerG1987b}. 
For example, with $B=3.454$, 
for $\lambda=1+1/(2B)$ 
we get a gain of  $\approx 6.1974568\cdot 10^{-7}$
obtained for $\ell=12$, $\omega =2491/2250$
and $\nu= 0.666856\dotsc$. 

In the proof of Theorem \ref{Thm-B-Cardinality} we also show  that
\eqref{B1-lowerbound} holds for 
\begin{equation}
\label{lambda-range}
\frac{1-\nu\Delta_{\ell,\omega}(\nu)}{1-\Delta_{\ell,\omega}(\nu)}
-
\frac{1-\Delta_{\ell,\omega}(\nu)}{2 B}
<
\lambda
\leq
\frac{1-\nu\Delta_{\ell,\omega}(\nu)}{1-\Delta_{\ell,\omega}(\nu)}
+
\frac{1-\Delta_{\ell,\omega}(\nu)}{2 B}
\end{equation}
thus extending, for $\Delta_{\ell,\omega}(\nu)\in (0,1)$ and 
$1/2<\nu<1-1/(2B)$, Cheer and Goldston's \cite{CheerG1987b}
result, which holds for $\lambda\in (1-1/(2B),1+1/(2B)]$,
to larger values of $\lambda$.

Finally,  in \eqref{dn-squared-lowerbound}
we get, for $\ell=12$, $\omega= 2491/2250$, $\nu\approx 0.666323\dotsc$
and $B=3.454$, the lower bound $1.00715710\dotsc$ which
improves the value $1.00698512\dotsc$ in eq.~(3.4) of Cheer and Goldston
\cite{CheerG1987b}; hence the constant $c_3(B)$ in
\eqref{dn-squared-lowerbound} can be chosen as $c_3(B)= 0.00017198\dotsc$.
%
%
%
%
Now we assume a suitable form of the $k$-tuple conjecture,
that is, equation \eqref{ktuple-conjecture} below, and, \emph{a
fortiori}, $B=1$ in Lemma \ref{BD-Thm2}.
The remark after the statement of Theorem \ref{Thm-A-Cardinality}
implies that Theorem \ref{Thm-B-Cardinality} and
\eqref{improved-lambda}-\eqref{lambda-range} still hold with the
condition on $\nu$ replaced by $\nu\in(0,1/2)$ and
$\Delta_{\ell,\omega}(\nu)$ replaced by 
$\widetilde{\Delta}_{\ell,\omega}(\nu)$ as defined in \eqref{Delta-tilde-def}.
In this case the numerical values for the key quantities are the following:
\begin{itemize}
\item $\lambda\geq 1.508146\dotsc$, for 
$\nu \approx 0.23435$, $\ell=10$ and $\omega\approx 39943/30000$,
to be compared with $\lambda=1.5$ in \cite{CheerG1987b};
%
%
\item $c_2(1,3/2) \approx 3.3185202\cdot10^{-5}$ for 
$\nu \approx 0.23435$, $\ell=10$ and $\omega\approx 39943/30000$,
to be compared with the value $0$ in \cite{CheerG1987b};
%
%
\item  \eqref{dn-squared-lowerbound} holds with the constant 
$1.09096653\dotsc$
obtained for $\nu \approx 0.22735$, $\ell=10$ and 
$\omega\approx 40027/30000$; this should
be compared with $1.08333333\dotsc$ in \cite{CheerG1987b}.
So we can choose $c_3(1)= 0.0076331\dotsc$.
%
\end{itemize}

The next result is about an arbitrary positive power of
$p_{i+1}-p_i$ whenever this distance is ``small''. We have the following 
\begin{Theorem}
\label{alpha-lower-bound}
Let $\epsilon>0$, $\alpha\geq 0$,
$X$ be a large parameter and $\lambda>1/2$.  Hence
\[
\sum_{\substack{p_{i}\leq X\\ p_{i+1}-p_{i}\leq \lambda\log X}}
(p_{i+1}-p_i)^\alpha 
\geq
(c_4(\lambda, \alpha)- \epsilon)
X (\log X)^{\alpha-1},
\]
where for $\alpha>0$  
\[
c_4(\lambda, \alpha)
=
\sup_{\ell \in \Z;\  \ell \geq 2} \sup_{\omega>1}
\left(
\frac{2\alpha\Delta_{\ell,\omega}(\lambda)}{(\alpha+1)B}
\right)^\alpha
\frac{\Delta_{\ell,\omega}(\lambda)}{\alpha+1},
\]
$\Delta_{\ell,\omega}(\lambda)$ is defined in \eqref{Delta-def},
$B$ is defined in Lemma \ref{BD-Thm2} below 
and  $c_4(\lambda, 0)= c_1(\lambda)$
is defined
in Theorem \ref{Thm-A-Cardinality}.
\end{Theorem}
 We remark that Theorem \ref{alpha-lower-bound} collapses to
Theorem \ref{Thm-A-Cardinality} for $\alpha=0$ and that
$c_4(\lambda, \alpha)$ is a decreasing function of $\alpha$. 
We also have a result concerning the set $\A_1(K)$ 
whenever $K$ is about $\log X$.
\begin{Theorem}
\label{A1-lower-bound}
Let $\epsilon>0$, $X$ be a sufficiently large parameter and
$0 < \lambda < 2 / B - \epsilon$
where $B$ is defined in Lemma \ref{BD-Thm2} below.  We have 
that
\[
\vert
\A_1(\lambda \log X)
\vert
\geq
\Bigl(
1- \lambda \frac{(B+\epsilon)}{2}
\Bigr)
\frac{X}{\log X}.
\]
\end{Theorem}
Recalling the result in Fouvry and Grupp \cite{FouvryG1986}, 
see also the remark after Lemma 
\ref{BD-Thm2}, we can set $B= 3.454$ 
and hence Theorem \ref{A1-lower-bound} holds for every  
$\lambda <  0.579038\dots$.
Assuming that the inequality in Lemma~\ref{BD-Thm2} holds with 
the best possible value $B = 1$ we obtain that Theorem \ref{A1-lower-bound} 
holds for every $\lambda <2$.
As a Corollary we have 
\begin{Corollary}
\label{B-lower-bound}
Let $\epsilon>0$, 
$X$ be a sufficiently large parameter and  $ \lambda>0$.  We have 
that
\[
\vert
\B(\lambda \log X)
\vert
\gg_{\lambda, \epsilon} 
X.
\]
\end{Corollary}
The implicit constant here is the same as in Theorem \ref{A1-lower-bound}
for $\lambda < 2 / B - \epsilon$ and it is $\epsilon$ 
otherwise.

Our last theorem is a conditional result on the cardinality 
of $\A_1(\lambda \log X)$ when  $\lambda$ is not in the range
described in Theorem \ref{A1-lower-bound}.
\begin{Theorem}
\label{A1-cond-lower-bound}
Let  $X$ be a sufficiently large parameter. 
For $\lambda \geq 2 / B$, where $B$ is defined in Lemma
\ref{BD-Thm2} below, we assume that
there exists a positive constant $c_5=c_5(\lambda)$
such that
\begin{equation}
\label{A1-cond-first-hyp}
\sum_{\substack{p_{i}\leq X\\ p_{i+1}-p_i > \lambda \log X }} 
(p_{i+1}-p_i) 
\geq 
c_5(\lambda) X.
\end{equation}
Assume further that there exists an absolute constant $c_6>0$
such that for every $\eta>\lambda$ we have
\begin{equation}
\label{A1-cond-second-hyp}
 \sum_{\substack{p_{i}\leq X\\  p_{i+1}-p_i > \eta \log X}} 
(p_{i+1}-p_i)
\leq
\frac{c_6}{\eta}
X.
\end{equation}
Then
\[
\vert
\A_1(\lambda \log X)
\vert
\gg_{\lambda,\eta} 
\frac{X}{\log X}.
\]
\end{Theorem}

We remark that Heath-Brown proved that the hypothesis
\eqref{A1-cond-second-hyp} holds under the assumption of the Riemann
Hypothesis and of a suitable form of the Montgomery pair-correlation
conjecture, see Corollary~1 of \cite{HeathBrown1982a}.
Theorem \ref{A1-cond-lower-bound} 
should be compared with 
Theorem 5 of Cheer and Goldston \cite{CheerG1987b}
in which our hypothesis \eqref{A1-cond-first-hyp} is replaced
by the stronger condition 
\[
\sum_{\substack{p_{i}\leq X\\ p_{i+1}-p_i > \lambda \log X }} 
(p_{i+1}-p_i - \lambda \log X) 
\gg_\lambda
X
\]
they used there. We finally remark that, using equation \eqref{max-constant}
below, we can also say that Theorem 5 of Cheer and Goldston \cite{CheerG1987b}
holds for $\lambda$ defined in \eqref{improved-lambda}.

\medskip
\paragraph{Acknowledgements.} 
We would like to thank Professors Heath-Brown and Maier 
for their insights and Professor  Perelli for an interesting discussion
about his joint papers with Salerno.

\section{Main Lemmas}
We will use two famous results
by Bombieri and Davenport.
\begin{Lemma}[Theorem 1 of Bombieri-Davenport \cite{BombieriD1966}]
\label{BD-Thm1}
Let  $1\leq T< (\log X)^C$ for some fixed positive $C$. Then, for any fixed 
positive $\epsilon$, we have
\[
\sum_{n=1}^T 
Z(X; 2n) 
> 
X
\sum_{n=1}^T \singseries(n)
-
\Bigl(
\frac14 +\epsilon
\Bigr) 
X\log X,
\]
where $Z(X; 2n)$ and $\singseries(n)$ are defined in 
\eqref{Z-def}-\eqref{singseries-def}.
\end{Lemma}
\begin{Lemma}[Theorem 2 of Bombieri-Davenport \cite{BombieriD1966}]
\label{BD-Thm2}
There exists a positive constant $B$ such that,
for any positive $\epsilon$ and 
for every positive integer $n$, we have
\[
Z(X; 2n) 
<
(B+\epsilon) \singseries(n) X,
\]
where $Z(X; 2n)$ and $\singseries(n)$ are defined in 
\eqref{Z-def}-\eqref{singseries-def}, provided $X$ is sufficiently large.
\end{Lemma}
Chen \cite{Chen1978} proved that $B=3.9171$
can be used in Lemma \ref{BD-Thm2}. Wu \cite{Wu2004} 
recently slightly improved this by 
proving that $B= 3.91045$ is admissible for every value of $n$.
For $n\leq \log^A X$, where $A>0$ is an arbitrary constant, 
the best result is $B= 3.454$ by Fouvry and Grupp \cite{FouvryG1986}.
Moreover we remark that similar results hold for $Z_1(X;2n)$,
defined in \eqref{Z-def}, since
$Z_1(X;2n) \leq Z(X;2n) \leq Z_1(X;2n) \log^2 X$ and hence, using
the inequalities
\[
\sum_{p\leq X/\log^4 X} 
\sum_{\substack{p'\leq X \\ p' -p =2n}} \log p \log p'
\leq
\pi
\Bigl(
\frac{X}{\log^4 X}
\Bigr) 
\log^2 X 
=
\odi{\frac{X}{\log^2 X}}
\]
and
\[
\sum_{X/\log^4 X<p \leq X} 
\sum_{\substack{p'\leq X \\ p' -p =2n}} \log p \log p'
>
(1+\odi{1})\log^2X 
\sum_{X/\log^4 X<p \leq X} 
\sum_{\substack{p'\leq X \\ p' -p =2n}} 1,
\]
we also obtain 
\begin{equation}
\label{Z1-Z-comparison}
\frac{Z(X;2n)}{\log^{2}X}
\leq
Z_1(X;2n) 
< 
(1 + \odi{1}) 
\frac{Z(X;2n)}{\log^{2}X}
+ 
\odi{\frac{X}{\log^2 X}}.
\end{equation}

Concerning the summation of the singular series of the twin-prime problem,
we will use the following result.
\begin{Lemma} [Friedlander-Goldston \cite{FriedlanderG1995}, eq.~(1.13)]
\label{singseries-averaged}
Let $X\geq 2$. We have
\[
\sum_{n\leq X} 
\singseries(n) 
=
X
+
\frac{1}{2} \log X + 
\Odi{(\log X)^{2/3}},
\]  
where $\singseries(n)$ is defined in \eqref{singseries-def}.
\end{Lemma}

\section{Proof of Theorem \ref{Thm-huge-costant} and Corollary
\ref{Corollary-log}}

We now define two different averages for primes in short intervals we
will need to prove Theorem \ref{Thm-huge-costant}. Let
\[
J_k(X,h) = \int_0^X (\psi(t+h) - \psi(t))^k \dx t
\quad
\textrm{and}
\quad
\widetilde{J}_k(X,h) = \sum_{m\leq X} (\psi(m+h) - \psi(m))^k
\]
be the Selberg integral and its discrete version.
In the proof of Theorem \ref{Thm-huge-costant} we will follow the line 
of Perelli and Salerno \cite{PerelliS1982,PerelliS1985} to connect the
moments over primes with the corresponding ones over 
integers.
To this end we now need several lemmas. 
We assume implicitly that $X$ is sufficiently large.

\begin{Lemma}[Gallagher \cite{Gallagher1976,Gallagher1981}]
\label{Gallagher}
Let $\epsilon>0$ and $1 \leq h\leq X$. Let further $k \geq 2$ be an integer. 
Then
\[  
\widetilde{J}_k(X,h) \leq 
\Bigl(
\P_k
\Bigl(
\frac{h}{\log X}
\Bigr)
+\epsilon
\Bigr)
X\log^k X.
\] 
where $\P_k(y)$ is defined in \eqref{P-def}.
\end{Lemma}
\begin{Proof}
The proof follows immediately inserting the following sieve estimate
of Klimov \cite{Klimov1958} in Gallagher's argument (see also Theorem
5.7 of \cite{HalberstamR1974})
\[
\sum_{m\leq X} 
\Lambda(m+h_1)
\Lambda(m+h_2)
\dotsm
\Lambda(m+h_r)
\leq 
(2^r r! + \epsilon) X \singseries(h_1,\dotsc,h_r)
\]
where $h_1,\dotsc,h_r$ are distinct integers such that  $0\leq h_i\leq h$
for every $i=1,\dotsc,r$,
\begin{equation}
\label{gen-singseries-def}
\singseries(h_1,\dotsc,h_r)
=
\prod_p 
\Bigl(
1-\frac{1}{p}
\Bigr)^{-r}
\Bigl(
1-\frac{\nu_p(h_1,\dotsc,h_r)}{p}
\Bigr)
\end{equation}
and $\nu_p(h_1,\dotsc,h_r)$ is the number of distinct residue classes
modulo $p$ the $h_i$, $i=1\dotsc,r$, occupy. 
\end{Proof}
The case $k=4$ of Lemma \ref{Gallagher} was recently proved 
by Goldston and Y{\i}ld{\i}r{\i}m, see 
eq.~(7.32) of  \cite{GoldstonY2007}.

We also remark that, assuming the $k$-tuple conjecture in the form
\begin{equation}
\label{ktuple-conjecture}
\sum_{m\leq X} 
\Lambda(m+h_1)
\Lambda(m+h_2)
\dotsm
\Lambda(m+h_r)
\sim
X \singseries(h_1,\dotsc,h_r)
\quad 
\textrm{as}
\ 
X\to +\infty
\end{equation}
where $h_1,\dotsc,h_r$ are distinct integers
and $\singseries(h_1,\dotsc,h_r)$
is defined in \eqref{gen-singseries-def},
Gallagher \cite{Gallagher1976} proved that 
 Lemma \ref{Gallagher} holds with the term
 $\P_k(h/\log X)$ replaced by $\widetilde{\P}_k (h/\log X)$ where
\begin{equation}
\label{P-tilde-def}
\widetilde{\P}_k (y)
=
\sum_{r=1}^{k} 
\Stirling{k}{r} y^r.
\end{equation}
Now we see two lemmas on the connections between 
the Selberg integral and its discrete version. 
\begin{Lemma}
\label{J-discreto-continuo1}
Let $h$ be an integer, $1\leq h\leq X$.  Let further $k \geq 2$ be an integer.  
Then
\[
J_k(X,h) = \widetilde{J}_k(X,h) + \Odi{\frac{h^k \log^k X}{\log^k (2h)}}.
\] 
\end{Lemma}
\begin{Proof}
It is clear that $\widetilde{J}_k(X,h) = 
\widetilde{J}_k(\lfloor X \rfloor,h) = J_k(\lfloor X \rfloor,h)$. 
By the Brun-Titchmarsh inequality we have $ J_k(X,h)= J_k(\lfloor X \rfloor,h)+ 
\Odi{h^k \log^k X(\log (2h))^{-k}}$ and the Lemma follows.
\end{Proof}
\begin{Lemma}
\label{J-discreto-continuo2}
Let $\epsilon>0$ and $h\in \R\setminus\Z$ with $1<h\leq X$.
Let further $k \geq 2$ be an integer.  
Then
\[
J_k(X,h) 
= 
\widetilde{J}_k(X,h) 
+ 
\Odip{k}{
X (\log X)^{k-1/2} 
\Bigl(
\P_{2k-2}
\Bigl(
\frac{2h}{\log X}
\Bigr)
+\epsilon
\Bigr)^{1/2}
}
+
\Odi{\frac{h^k \log^k X}{\log^k (2h)}},
\]
where  $\P_k(y)$ is defined in \eqref{P-def}.
\end{Lemma}
\begin{Proof}
By the Brun-Titchmarsh inequality and
letting $h=\lfloor h \rfloor+\beta$ and $t=\lfloor t \rfloor+\tau$ we have
\[
\begin{split}
J_k(X,h) 
&= 
\int_0^X 
(
\psi(\lfloor t \rfloor+\lfloor h \rfloor+\beta+\tau) 
- 
\psi(\lfloor t \rfloor+\tau)
)^k 
\dx t
\\
&=
\sum_{m\leq X} 
\int _0^1
(
\psi(m+\lfloor h \rfloor+\beta+\tau) 
-
\psi(m)
)^k 
\dx \tau
+ 
\Odi{\frac{h^k \log^k X}{\log^k (2h)}}
\\
&=
\sum_{m\leq X} 
\int _0^{1-\beta}
(\psi(m+\lfloor h \rfloor) - \psi(m))^k \dx \tau
+
\sum_{m\leq X} 
\int _{1-\beta}^1
(\psi(m+\lfloor h \rfloor+1) - \psi(m))^k \dx \tau 
\\
& 
\hskip2cm
+ \Odi{\frac{h^k \log^k X}{\log^k (2h)}} 
\\
&=
(1-\beta)  \widetilde{J}_k(X,\lfloor h \rfloor) + 
\beta\  \widetilde{J}_k(X,\lfloor h \rfloor+1)+ 
\Odi{\frac{h^k \log^k X}{\log^k (2h)}}.
\end{split}
\]
We have
\[
  \widetilde{J}_k(X, \lfloor h \rfloor+ 1)
  -
  \widetilde{J}_k(X, \lfloor h \rfloor) 
  =
  \sum_{m \le X}
    \biggl(
      (\psi(m + \lfloor h \rfloor + 1) - \psi(m))^k
      -
      (\psi(m + \lfloor h \rfloor) - \psi(m))^k
    \biggr).
\]
Let $a = \psi(m + \lfloor h \rfloor + 1) - \psi(m)$ and
$b = \psi(m + \lfloor h \rfloor) - \psi(m)$: By the Mean Value Theorem
we have $a^k - b^k \le k a^{k-1} (a - b)$, since $b \le a$.
By the Cauchy-Schwarz inequality, the Prime Number Theorem and the identity
$\psi(m + \lfloor h \rfloor + 1) - \psi(m + \lfloor h \rfloor) 
= \Lambda(m + \lfloor h \rfloor + 1)$, we get
\[
\begin{split}
  \widetilde{J}_k(X, \lfloor h \rfloor+ 1)
  -
  \widetilde{J}_k(X, \lfloor h \rfloor) 
&  
\ll_k
  \Bigl(  J_{2k-2}(X,\lfloor h \rfloor+1)\Bigr)^{1/2}
  \Bigl( \sum_{m \le X} \Lambda(m + \lfloor h \rfloor + 1)^2 \Bigr)^{1/2}
\\
&
  \ll_k
  (X \log X)^{1/2}
\Bigl[    
  X \log^{2k-2} X 
  \Bigl(
\P_{2k-2}
\Bigl(
\frac{\lfloor h \rfloor+1}{\log X}
\Bigr)
+\epsilon
\Bigr)
\Bigr]^{1/2}\\
&\ll
X (\log X)^{k-1/2} 
 \Bigl(
\P_{2k-2}
\Bigl(
\frac{2h}{\log X}
\Bigr)
+\epsilon
\Bigr)^{1/2}.
\end{split}
\]
Lemma \ref{J-discreto-continuo2} now follows because
$\widetilde{J}_k(X,h) = \widetilde{J}_k(X,\lfloor h \rfloor)$.
\end{Proof}

Let now $u$ be a positive real number and 
\begin{equation}
\label{def-psik}
\psi_k (X, u) = 
\sum_{\substack{{m_1,\dotsc, m_k}\\ 
\min(m_i) \leq X \\ \max(m_i)-\min(m_i) \leq u}}
\Lambda(m_1) \dotsm \Lambda(m_k).
\end{equation}
This function can be easily connected with the Selberg integral.
\begin{Lemma}[Perelli-Salerno \cite{PerelliS1982,PerelliS1985}]
\label{media-kuple}
Let $1 \leq h \leq X $, and $k \geq 2$ be an integer. 
Then
\[
J_{k}(X,h) 
=
\int_0^h 
\psi_k (X, u) \dx u
+\Odi{\frac{h^{k+1} \log^k X}{\log^k (2h)}}.
\]
\end{Lemma}
\begin{Proof}
Let $N=\max(m_i)$ and $n=\min(m_i)$ in \eqref{def-psik}.
Expanding the $k$-th power in $J_{k}(X,h)$, we have
\[
\begin{split}
J_{k}(X,h) 
&= 
\sum_{\substack{{m_1,\dotsc, m_k}\\ n \leq X \\ N-n \leq h}}
\Lambda(m_1) \dotsm \Lambda(m_k) (h-N+n)
\\
&+
\sum_{\substack{{m_1,\dotsc, m_k}\\  X < n \leq X+h \\ N-n \leq h}}
\Lambda(m_1) \dotsm \Lambda(m_k) (X+h-N)
=
\Sigma_1+\Sigma_2,
\end{split}
\]
say.
By the partial summation formula we 
immediately get
\(
\Sigma_1 = 
\int_0^h 
\psi_k (x, u) \dx u
\)
and, by the Brun-Titchmarsh inequality, we have
\(
\Sigma_2 \ll  h^{k+1} \log^k X (\log (2h))^{-k}.
\)
Lemma \ref{media-kuple} follows.
\end{Proof}
The next lemma gives an upper bound for $\psi_k(X,h)$ in
terms of the discrete Selberg integral.
\begin{Lemma}
\label{monotonicity}
Let $\epsilon>0$, $\omega>1$ and $1 \leq h\leq X$.
Let further $k \geq 2$ be an integer. 
Then
\[
\psi_k (X, h)
\leq 
\frac
{\widetilde{J}_k(X,\omega\, h)}
{(\omega-1)h}
+
\Odi{\frac{\omega^{k+1} h^{k} \log^k X}{(\omega-1)\log^k (2\omega\, h)}}
+
\Odiprimep{k}{
\frac{X\log^{k-1/2} X}{(\omega-1)h} 
 \Bigl(
\P_{2k-2}
\Bigl(
\frac{2\omega\, h}{\log X}
\Bigr)
+\epsilon
\Bigr)^{1/2}
}
\]
where $\mathcal{O'}$ means that this error term is present 
only if $h\not \in \Z$ and $\P_k(y)$ is defined in \eqref{P-def}.
\end{Lemma}
\begin{Proof}
Since $\psi_k (X, u)$ is a positive and increasing function of $u$,
it is easy to see that
\[
\psi_k (X, h)
\leq 
\frac{1}{(\omega-1)h}
\int_h^{\omega\, h} \psi_k (X, u) \dx u
\leq 
\frac{1}{(\omega-1)h}
\int_0^{\omega\, h} \psi_k (X, u) \dx u,
\]
where $\omega>1$ is a constant.
By Lemmas \ref{J-discreto-continuo1}, \ref{J-discreto-continuo2}
and \ref{media-kuple}, we have
\[
\begin{split}
\int_0^{\omega\, h} \psi_k (X, u) \dx u
\ = \ 
&
\widetilde{J}_k(X,\omega\, h)
+
\Odi{\frac{\omega^{k+1} h^{k+1} \log^k X}{\log^k (2\omega\, h)}}
\\
&
+
\Odiprimep{k}{
X\log^{k-1/2} X 
\Bigl(
\P_{2k-2}
\Bigl(
\frac{2\omega\, h}{\log X}
\Bigr)
+\epsilon
\Bigr)^{1/2}
}
\end{split}
\]
and hence Lemma \ref{monotonicity} follows.
\end{Proof}
The following last lemma is a lower bound for $\psi_{k+1}(X,h)$ in
terms of a weighted form of the discrete Selberg integral.
\begin{Lemma}[Perelli-Salerno \cite{PerelliS1982,PerelliS1985}]
\label{average-over-primes}
Let $1 \leq h \leq X$, and $k \geq 2$ be an integer. 
Then
\[
\sum_{m\leq X}
\Lambda(m)
(\psi(m+h)-\psi(m))^k
\leq
\psi_{k+1}(X,h).
\]
\end{Lemma}
\begin{Proof}
First of all we remark that
\begin{equation}
\label{equival}
\sum_{m\leq X}
\Lambda(m)
(\psi(m+h)-\psi(m))^k
=
\sum_{m\leq X}
\Lambda(m)
\sum_{\substack{{m_1,\dotsc, m_k}\\ m < m_i \leq m+h}}
\Lambda(m_1) \dotsm \Lambda(m_k).
\end{equation}
Recalling now that $N=\max(m_i)$ and $n=\min(m_i)$ in \eqref{def-psik},
we trivially have
\begin{equation}
\label{lower-bound}
\psi_{k+1}(X,h)  
\geq 
\sum_{n\leq X} \Lambda(n)
\sum_{\substack{m_2, \dotsc, m_{k+1}\\ m_i \neq n \\ 0 < N-n \leq h}}
\Lambda(m_2) \dotsm \Lambda(m_{k+1}) 
=
\sum_{n\leq X}
\Lambda(n)
\sum_{\substack{{m_1,\dotsc, m_k}\\ n < m_i \leq n+h}}
\Lambda(m_1) \dotsm \Lambda(m_k)
\end{equation}
and Lemma \ref{average-over-primes} follows immediately 
combining \eqref{equival} and \eqref{lower-bound}.
\end{Proof}

\smallskip
Now we are ready to prove Theorem \ref{Thm-huge-costant} and Corollary
\ref{Corollary-log}.

By Lemmas \ref{average-over-primes}, \ref{monotonicity} and \ref{Gallagher} 
the upper bound in the statement of Theorem \ref{Thm-huge-costant}
holds for the function 
$ (\log X)^{-1} \sum_{m\leq X} \Lambda(m) (\psi(m+h)-\psi(m))^k$.
It is easy to see that the contribution of $m=p^\alpha$ with $\alpha>1$
in the previous sum is negligible. Theorem \ref{Thm-huge-costant}  hence
follows by the partial summation formula since $f(X) \leq  h\leq X^{1-\epsilon}$,
where $f(X) \to +\infty$ arbitrarily slowly as $X\to +\infty$.
The first part of Corollary \ref{Corollary-log} can be obtained 
inserting $h=\lambda \log X$ in 
Theorem \ref{Thm-huge-costant}
while the second part follows from the first one using the partial 
summation formula.

\section{Proof of Theorem \ref {Thm-A-Cardinality}}
Let $X$ be a large parameter and $2\leq K\leq X$.
Our goal here is to elementarily prove a lower bound of 
the correct order of magnitude for  the cardinality of the set 
$\A(K)$ defined in \eqref{A(K)-def} when $K$ is about $\log X$.
Let $\ell \geq 2$ be an integer. By 
the H\"older inequality we obtain   
\[
\sum_{p\in \A(K)} 
(\pi(p+K)-\pi(p)) 
\leq
\Bigl(\sum_{p\in \A(K)} 
(\pi(p+K)-\pi(p))^\ell\Bigr)^{1/\ell}
\vert \A(K) \vert^{(\ell-1)/\ell}
\]
and hence
\begin{equation} 
\label{B-lower1}
\vert \A(K) \vert  
\geq 
\Bigl[
\sum_{p\in \A(K)} (\pi(p+K)-\pi(p)) 
\Bigr]^{\ell/(\ell-1)}
\Bigl(\sum_{p\in \A(K)} 
(\pi(p+K)-\pi(p))^\ell\Bigr)^{-1/(\ell-1)}.
\end{equation}
Now we proceed to estimate the term in the square brackets.
It is easy to see that 
\begin{equation} 
\label{A-Z1-link}
\sum_{p\in \A(K)} 
(\pi(p+K)-\pi(p))
=
\sum_{p \leq X} 
(\pi(p+K)-\pi(p))  
=
\sum_{n=1}^{K/2} 
Z_1(X; 2n),
\end{equation}
where $Z_1(X; 2n)$ is defined in \eqref{Z-def}, because
$\pi(p + K) = \pi(p)$ by definition if $p \notin \A(K)$.
From Lemmas \ref{BD-Thm1} and \ref{singseries-averaged}
and \eqref{Z1-Z-comparison}, 
we get
\[
\sum_{n=1}^{K/2} Z_1(X; 2n) 
> 
\frac{KX}{2\log^2 X} 
+
\frac{X\log (K/2)}{2\log^2 X}
-
\Bigl(
\frac14 +\epsilon
\Bigr) 
\frac{X}{\log X}
+
\Odi{X \frac{(\log (K/2))^{2/3}}{\log^2 X}}.
\]
Hence, if $K=\lambda \log X$, with $\lambda>1/2$
and $X$ is sufficiently large, the previous remark implies
\begin{equation}
\label{sq-brackets-estim}
\sum_{p\in \A(\lambda\log X)} 
(\pi(p+ \lambda \log X)-\pi(p))
> 
\Bigl(
\frac{\lambda}{2} -\frac14 -\epsilon/2
\Bigr) 
\frac{X}{\log X}.
\end{equation}
To estimate the second term in \eqref{B-lower1}, we use 
the second part of Corollary \ref{Corollary-log} with $k=\ell$. 
We immediately get
\begin{equation}
\label{square-estim}
\sum_{p \in \A(\lambda\log X)}
(\pi(p+\lambda\log X) -\pi(p))^\ell
=
\sum_{p \leq X}
(\pi(p+\lambda\log X) -\pi(p))^\ell
\ll_{\ell, \lambda, \omega, \epsilon}
\frac{X}{\log  X},
\end{equation}
where the implicit constant, using Lemmas \ref{Gallagher}, \ref{J-discreto-continuo2}
and \ref{monotonicity}, is  $\RR_{\ \ell+1,\omega}(\lambda)+\epsilon$ with
$\RR_{\ \ell,\omega}(\lambda)$ is defined in \eqref{R-def}.
Inserting now \eqref{sq-brackets-estim}-\eqref{square-estim} in \eqref{B-lower1} we obtain Theorem \ref {Thm-A-Cardinality} with the implicit constant 
equal to
$c(\ell, \lambda, \omega, \epsilon)
= 
(\lambda/2 -1/4)^{\ell/(\ell-1)}
\RR_{\ \ell+1,\omega}(\lambda)^{-1/(\ell-1)}
-\epsilon$
for every $\omega>1$ and every integer $\ell \geq 2$. 
Hence the best possible constant is
$\sup_{\ell \in \Z;\  \ell \geq 2}  
\sup_{\omega>1}
c(\ell, \lambda,\omega,\epsilon)$.

Using the $k$-tuple conjecture, we can use \eqref{P-tilde-def}
and, moreover, \eqref{sq-brackets-estim} holds
with $\lambda/2$ instead of $\lambda/2-1/4$.
Hence the implicit constant in this theorem is equal to 
$\widetilde{c}(\ell, \lambda, \omega,\epsilon)
= 
(\lambda/2)^{\ell/(\ell-1)}
\widetilde{\RR}_{\ \ell+1,\omega}(\lambda)^{-1/(\ell-1)}
-\epsilon$,
where $\widetilde{\RR}_{\ \ell,\omega}(\lambda)$ is
defined in \eqref{R-tilde-def}.

\section{Proof of Theorem \ref{Thm-B-Cardinality}}

Let $1/2 < \mu < \lambda$ be real numbers.
By the Prime Number Theorem 
we get $p_{\pi(X)+1}= X(1+\odi{1})$ and hence,
following the line of Cheer and Goldston's proof of Theorem 1 
in \cite{CheerG1987b},
we get that
\begin{equation}
\label{CG-startingpoint}
\begin{split}
1 + \odi{1}
&
=
\frac{p_{\pi(X)+1}-2}{X} 
=
\frac{1}{X}
\sum_{p_{i}\leq X}
(p_{i+1}-p_i) 
\\
&= 
\mu 
+ 
\frac{1}{X}
\Bigl(
- 
S_X(\mu) 
+ 
\sum_{\substack{
p_{i}\leq X\\ 
\mu  \log X< p_{i+1}-p_i \leq \lambda \log X }} 
(p_{i+1}-p_i - \mu \log X) 
\\
& \hskip2cm
+
\sum_{\substack{
p_{i}\leq X\\ 
p_{i+1}-p_i > \lambda \log X }} 
(p_{i+1}-p_i - \mu \log X)
\Bigr) + o(1),
\end{split}
\end{equation}
where 
\[
S_X(\mu) 
= 
\sum_{\substack{
p_{i}\leq X\\ 
p_{i+1}-p_i \leq \mu \log X }} 
( \mu \log X - p_{i+1} + p_i ).
\]
Choose $\nu \in (1/2, \mu)$: by Theorem \ref{Thm-A-Cardinality} we
have
\begin{equation}
\label{S-lowerbound}
S_X(\mu) 
\geq 
\sum_{\substack{
p_{i}\leq X\\ 
p_{i+1}-p_i \leq \nu \log X }} 
( \mu \log X - p_{i+1} + p_i )
>
(\mu-\nu) \log X
\sum_{\substack{
p_{i}\leq X\\ 
p_{i+1}-p_i \leq \nu \log X }} 
1
>
(\mu-\nu) 
\Delta_{\ell,\omega}(\nu) X,
\end{equation}
where $\Delta_{\ell,\omega}(\nu)$ is defined in \eqref{Delta-def}.
Inserting \eqref{S-lowerbound} in \eqref{CG-startingpoint} and
arguing again as on page 475 of Cheer and Goldston \cite{CheerG1987b}, 
we get 
\begin{equation}
\label{better}
1
- \mu 
+ (\mu-\nu) \Delta_{\ell,\omega}(\nu) 
- B \frac{(\lambda-\mu)^2}{2}
- \epsilon
\leq 
\frac{1}{X}
\sum_{\substack{
p_{i}\leq X\\ 
p_{i+1}-p_i > \lambda \log X }} 
(p_{i+1}-p_i - \mu \log X)
\end{equation}
and we gain the summand $(\mu-\nu) \Delta_{\ell,\omega}(\nu)$ 
comparing this equation with (3.8) in \cite{CheerG1987b}.
Now we have to optimize the LHS of \eqref{better}.
If we consider $\lambda$ and $\nu$ fixed then the maximum is attained for
$\mu = \lambda - (1 - \Delta_{\ell,\omega}(\nu)) / B$ and this, arguing as
on page 476 of Cheer and Goldston \cite{CheerG1987b}, leads to
\[
  \lambda
  <
  \Bigl( 1 - \Delta_{\ell,\omega}(\nu) \Bigr)^{-1}
  \Bigl(
    1
    - \nu \Delta_{\ell,\omega}(\nu)
    + \frac{(1 - \Delta_{\ell,\omega}(\nu))^2}{2 B}
  \Bigr)
  -
  \epsilon
\]
for every $\omega>1$ and every integer $\ell \geq 2$.
Developing the right hand side in powers of $\Delta_{\ell,\omega}$,
which is small, we see that we have
\[
  \lambda
  <
  1
  +
  \frac1{2 B}
  +
  \Delta_{\ell,\omega}(\nu)
  \Bigl(1 - \frac1{2 B} - \nu \Bigr)
  +
  \Odi{\Delta_{\ell,\omega}(\nu)^2}.
\]
For $1/2 < \nu< 1- 1 / (2 B)$, our net gain over the result of Cheer
and Goldston \cite{CheerG1987b} is essentially in the third summand
above, which we may maximize over $\omega$, $\ell$ and $\nu$.

Now we estimate $c_2(B,\lambda)$.
First of all we recall that Gallagher, see the remark at
the bottom of p.~87
of \cite{HeathBrown1982a}, 
proved that
\begin{equation}
\label{gallagher-estimate}
c_2(B,\lambda)
\geq
1-\lambda,
\end{equation}
for any $\lambda \in (0,1)$.
Let now 
$\lambda
<
\Bigl( 1 - \Delta_{\ell,\omega}(\nu) \Bigr)^{-1}
  \Bigl(
    1
    - \nu \Delta_{\ell,\omega}(\nu)
    + (1 - \Delta_{\ell,\omega}(\nu))^2/(2 B)
  \Bigr)
  -
  \epsilon
$.
Starting again from \eqref{better}, arguing as at the bottom of
p.~476 of \cite{CheerG1987b} and letting $\tau \geq \lambda$, we get
\begin{equation}
\label{max-constant}
\begin{split}
\frac{1}{X}
&
\sum_{\substack{
p_{i}\leq X\\ 
p_{i+1}-p_i > \lambda \log X }} 
(p_{i+1}-p_i - \lambda \log X)
\geq
\\
&
\Bigl[
1
- \mu 
+ (\mu-\nu) \Delta_{\ell,\omega}(\nu) 
- B \frac{(\lambda-\mu)^2}{2}
- \frac{B}{2} (\lambda-\mu)(\tau-\lambda)
\Bigr]
\Bigl(
1
+
\frac{\lambda-\mu}{\tau-\lambda}
\Bigr)^{-1}
  -
  \epsilon.
\end{split}
\end{equation}
Now, given $\lambda$ and $\nu$, we would like to maximize this term with
respect to $\tau$ and $\mu$.  
Using the substitution $u_1=\lambda-\mu>0$ and $u_2=(\tau-\lambda)/u_1\geq0$,
the RHS of \eqref{max-constant} becomes
\begin{equation}
\label{RHS}
\Bigl[
1 + u_1 - \lambda
-
\frac{B}{2} (1+u_2)u_1^2
+
(\lambda-u_1-\nu) \Delta_{\ell,\omega}(\nu) 
\Bigr]
\frac{u_2}{u_2+1}.
\end{equation}
For $\Delta_{\ell,\omega}(\nu)\in (0,1)$
and $\lambda \in (\lambda_1,\lambda_2]$,
where
\begin{equation}
\label{lambda-interval}
\lambda_1
=
\frac{1-\nu\Delta_{\ell,\omega}(\nu)}{1-\Delta_{\ell,\omega}(\nu)}
-
\frac{1-\Delta_{\ell,\omega}(\nu)}{2 B}
\quad
\textrm{and}
\quad
\lambda_2
=
\frac{1-\nu\Delta_{\ell,\omega}(\nu)}{1-\Delta_{\ell,\omega}(\nu)}
+
\frac{1-\Delta_{\ell,\omega}(\nu)}{2 B}
,
\end{equation}
equation \eqref{RHS} is maximized by
\[
  u_1
  =
  \frac{1-\Delta_{\ell,\omega}(\nu)}{B (1 + u_2)}
\]
and
\[
u_2= 
\frac
{(1-\Delta_{\ell,\omega}(\nu))^2
-
2 B (\lambda-1-(\lambda-\nu)\Delta_{\ell,\omega}(\nu))
}
{(1-\Delta_{\ell,\omega}(\nu))^2
+
2 B (\lambda-1-(\lambda-\nu)\Delta_{\ell,\omega}(\nu))
}.
\]
For these values of $u_1,u_2$, the maximal  RHS 
of \eqref{max-constant} is
\begin{equation}
\label{best-constant}
\frac{B}{2(1-\Delta_{\ell,\omega}(\nu))^2}
\Bigl(
1-\lambda
+
\frac{(1-\Delta_{\ell,\omega}(\nu))^2}{2 B}
+
(\lambda-\nu)\Delta_{\ell,\omega}(\nu)
\Bigr)^2
\end{equation}
which is larger than \eqref{gallagher-estimate}
for every $\lambda\in(\lambda_1,\lambda_2]$,
where $\lambda_1,\lambda_2$ are defined in \eqref{lambda-interval}.

Moreover, since $\Delta_{\ell,\omega}(\nu)\in (0,1)$ for
$1/2<\nu<1-1/(2B)$, we have that equation \eqref{lambda-interval} 
extends both the
width of $\lambda$-interval and the lower bound of eq.~(3.3) of
Cheer and Goldston \cite{CheerG1987b} to larger values.
Comparing again with eq.~(3.3) of Cheer and Goldston \cite{CheerG1987b},
for $1/2<\nu<\lambda-1/B$ our equation \eqref{best-constant} gives
a larger value for the final constant. 
This completes the proof of \eqref{B1-lowerbound}.

To prove \eqref{dn-squared-lowerbound} we follow again the line
of the proof of Theorem 3 of Cheer and Goldston \cite{CheerG1987b};
the only difference is paying attention to split the integration
interval in eq.~(3.10) there into the subintervals
$[0,1-1/(2B)]$, $[1-1/(2B), \lambda_1]$, $[\lambda_1, \lambda_2]$ and
where $\lambda_1$ and $\lambda_2$ are defined in \eqref{lambda-interval}.
We remark that $\lambda_2 \ge 1 + 1 / (2 B)$ since
$1/2< \nu < 1 - 1 / (2 B)$.
We have
\begin{align*}
  \sum_{p_i \le X} (p_{i+1}- p_{i})^2
  &\ge
  2 (1 - \epsilon)
  X \log X \int_0^{\lambda_2} 
  c_1(B, \lambda) \, \dx \lambda \\
  &=
  (1 - \epsilon) X \log X
  \Bigl\{
    2 \int_0^{1 - 1/(2B)}
    +
    2 \int_{1 - 1/(2B)}^{\lambda_1}
    +
    2 \int_{\lambda_1}^{\lambda_2}
    c_1(B, \lambda)
    \dx \lambda
  \Bigr\} \\
  &=
  (1 - \epsilon) X \log X (I_1 + I_2 + I_3),
\end{align*}
say, where in the first integral we use Gallagher's estimate
\eqref{gallagher-estimate}, in the second Cheer and Goldston's
equation~(3.3) from \cite{CheerG1987b}, and in the third our lower
bound~\eqref{best-constant}.
A fairly tedious computation reveals that
\begin{align*}
  I_1
  &\geq	
  1 - \frac1{4 B^2}, \\
  I_2
  &\geq
  \frac1{3 B^2}
  -
  \frac B3
  \Bigl( \frac{\Delta_{\ell,\omega}(\nu) (\nu - 1)}
              {1 - \Delta_{\ell,\omega}(\nu)}
         +
         \frac{2 - \Delta_{\ell,\omega}(\nu)}{2 B}
  \Bigr)^3, \\
  I_3
  &\geq
  \frac{(1 - \Delta_{\ell,\omega}(\nu))^3}{3 B^2}.
\end{align*}
Notice that $\Delta_{\ell,\omega}(\nu) = 0$ yields exactly
equation~(3.4) in \cite{CheerG1987b}.
Developing again in powers of $\Delta_{\ell,\omega}$, we get
\[
  I_1 + I_2 + I_3
  \ge
  1
  +
  \frac1{12 B^2}
  +
  \frac{\Delta_{\ell,\omega}(\nu)}B
  \Bigl( 1 - \frac1{2 B} - \nu  \Bigr)
  +
  \Odi{\Delta_{\ell,\omega}(\nu)^2}.
\]
Comparing this with equation~(3.4) of \cite{CheerG1987b},
we see that our gain comes from the third term above which is 
positive for $\nu<1-1/(2B)$.
This completes the proof of Theorem \ref{Thm-B-Cardinality}.

\section{Proof of Theorem \ref{alpha-lower-bound}}

Let $X$ be a large parameter and
$0<\eta<\lambda$. First we remark that the case $\alpha = 0$ 
corresponds to Theorem \ref{Thm-A-Cardinality}.
So from now on we can assume $\alpha>0$. We have
\[
\begin{split}
&
\sum_{\substack{p_{i}\leq X\\ p_{i+1}-p_{i}\leq \lambda\log X}}
(p_{i+1}-p_{i})^\alpha 
\geq
\sum_{\substack{p_{i}\leq X\\ \eta \log X < p_{i+1}-p_{i}\leq \lambda\log X}}
(p_{i+1}-p_{i})^\alpha 
 \\
&
> 
(\eta \log X)^\alpha 
\sum_{\substack{p_{i}\leq X\\ \eta \log X < p_{i+1}-p_{i}\leq \lambda\log X}} 
1
= 
(\eta \log X)^\alpha 
\Bigl(
\sum_{\substack{p_{i}\leq X\\ p_{i+1}-p_{i}\leq \lambda\log X}}
1  
- 
\sum_{\substack{p_{i}\leq X\\ p_{i+1}-p_{i}\leq \eta\log X}} 
1 
\Bigr).
\end{split}
\]
For $\lambda >1/2$ we can apply 
Theorem \ref{Thm-A-Cardinality} thus getting
\begin{equation}
\label{S-lower}
\sum_{\substack{p_{i}\leq X\\ p_{i+1}-p_{i}\leq \lambda\log X}}
(p_{i+1}-p_{i})^\alpha 
> 
(\eta \log X)^\alpha
\Bigl(
 \frac{X}{\log X} 
 ( \Delta_{\ell,\omega}(\lambda) -\epsilon)
 - 
\sum_{\substack{p_{i}\leq X\\ p_{i+1}-p_{i}\leq \eta\log X}} 
1 
\Bigr),
\end{equation}
where $ \Delta_{\ell,\omega}(\lambda)$ is defined in  \eqref{Delta-def}.
Using now equation \eqref{Z1-Z-comparison} and 
Lemmas \ref{BD-Thm2}-\ref{singseries-averaged}, we get
\begin{equation}
\label{bound-Z1}
\sum_{\substack{p_{i}\leq X\\ p_{i+1}-p_{i}\leq \eta\log X}} 
1 
\leq 
\sum_{n\leq (\eta/2) \log X} 
Z_1(X,2n) 
< 
\frac{\eta}{2} (B+\epsilon) 
\frac{X}{\log X}
\end{equation}
and hence, 
by \eqref{S-lower}-\eqref{bound-Z1}, we obtain
\[ 
\sum_{\substack{p_{i}\leq X\\ p_{i+1}-p_{i}\leq \lambda\log X}}
(p_{i+1}-p_{i})^\alpha 
> 
X\eta^\alpha (\log X)^{\alpha-1}
\Bigl(
\Delta_{\ell,\omega}(\lambda) 
 - 
\frac{\eta}{2} (B+\epsilon) 
\Bigr).
\]
Choosing the optimal value
$
\eta
= 
2\alpha
\Delta_{\ell,\omega}(\lambda)/( (\alpha+1) B)
$
for $\alpha>0$, Theorem \ref{alpha-lower-bound} follows.

\section{Proof of Theorem \ref{A1-lower-bound} and Corollary \ref{B-lower-bound}}

Let $X$ be a large parameter and $2\leq K\leq X$.
Arguing as in \eqref{A-Z1-link}  we obtain
\[
\vert 
\A(K)
\vert
=
\sum_{p\in \A(K)} 1
\leq
\sum_{p\in \A(K)} 
(\pi(p+K)-\pi(p)) 
=
\sum_{n=1}^{K/2} 
Z_1(X; 2n)
\]
and hence, by Lemmas  \ref{BD-Thm2}  
and \ref{singseries-averaged} we get 
\[
\vert 
\A(K)
\vert
\leq 
(B+\epsilon)
\frac{X}{\log^2 X}
\sum_{n=1}^{K/2} 
\singseries(n)
<
(B+\epsilon)
\frac{KX}{2\log^2 X}.
\]
So we have that 
$
\vert \A_1(K) \vert
>
\pi(X) - 
((B+\epsilon)/2)
X K
 (\log X)^{-2}
$
and hence, letting $\lambda>0 $ and $K=\lambda \log X$, 
we immediately get
\[
\vert \A_1(\lambda \log X) \vert
\geq 
\pi(X)
- 
(B+\epsilon)
\lambda
\frac{X}{2\log X}
>
  \Bigl( 1 - (B+\epsilon)\frac{\lambda}{2} \Bigr)
\frac{X}{\log X}
\]
for every sufficiently large $X$.
In the last inequality we have used 
$\pi(X)> X  (\log X)^{-1}$  for every $X\geq 17$ proved in 
Corollary 1 of Rosser-Schoenfeld \cite{RosserS1962}.
Choosing $\lambda < 2 / B - \epsilon$, Theorem \ref{A1-lower-bound}
follows at once. 

The proof of Corollary \ref{B-lower-bound} runs as follows.
The starting point is the trivial relation
\[
\vert
\B(K)
\vert
>
\sum_{\substack{p_{i}\leq X\\ p_{i}-p_{i-1}> K}} 
K 
+
\sum_{\substack{p_{i}\leq X\\ p_{i}-p_{i-1}\leq K}} 
(p_{i}-p_{i-1})
=
\sum_{\substack{p_{i}\leq X\\ p_{i+1}-p_{i}> K}} 
K 
+
\sum_{\substack{p_{i}\leq X\\ p_{i+1}-p_{i} \leq K}} 
(p_{i+1}-p_{i})
+
\Odi{K} .
\]
Letting $K=\lambda \log X$ and using Theorem \ref{A1-lower-bound}, 
we immediately get
\[
\vert
\B(\lambda \log X )
\vert
>
\lambda \log X 
\sum_{\substack{p_{i}\leq X\\ p_{i+1}-p_{i}> \lambda \log X}} 
1 
+
\Odip{\lambda}{\log X}
\gg_{\lambda,\epsilon}
X
\]
for every $0 < \lambda < 2 / B - \epsilon$.
Corollary \ref{B-lower-bound} now follows by remarking 
that 
$\vert
\B(\lambda \log X )
\vert
$
is an increasing function of $\lambda$.

\section{Proof of Theorem \ref{A1-cond-lower-bound}}

Let $\lambda \geq 2/B$ and $\eta>\lambda$ be a real number.
It is easy to see that
\begin{equation}
\label{starting-point}
\begin{split}
\sum_{\substack{p_{i}\leq X\\ p_{i+1}-p_i > \lambda \log X}} 
1 
&
\geq
\frac{1}{\eta \log X} 
\sum_{\substack{p_{i}\leq X\\ \lambda \log X < p_{i+1}-p_i \leq \eta \log X}} 
(p_{i+1}-p_i)
\\
&
=
\frac{1}{\eta \log X} 
\Bigl( 
\sum_{\substack{p_{i}\leq X\\  p_{i+1}-p_i > \lambda \log X }} 
(p_{i+1}-p_i)
- 
\sum_{\substack{p_{i}\leq X\\  p_{i+1}-p_i > \eta \log X }} 
(p_{i+1}-p_i)
\Bigr).
\end{split}
\end{equation}
Inserting the hypotheses \eqref{A1-cond-first-hyp} and  
\eqref{A1-cond-second-hyp} into \eqref{starting-point} we
immediately get
\[
\sum_{\substack{p_{i}\leq X\\ p_{i+1}-p_i > \lambda \log X}} 
1 
> 
\frac{X}{\eta\log X} 
\Bigl(
c_5(\lambda) -\frac{c_6}{\eta}
\Bigr).
\]
Theorem \ref{A1-cond-lower-bound} now follows since, for any fixed
$\lambda$, $\eta$ can be chosen greater than 
$\max(\lambda;c_6/c_5(\lambda))$.

\bibliographystyle{plain}

\bigskip
\noindent
\begin{tabular}{l@{\hskip 3mm}l@{\hskip 3mm}l}
D.~Bazzanella & A.~Languasco & A.~Zaccagnini \\
Politecnico di Torino & Universit\`a di Padova & Universit\`a di Parma \\
Dipartimento di Matematica & Dipartimento di Matematica 
& Dipartimento di Matematica \\
Corso Duca degli Abruzzi 24 & Pura e Applicata
& Parco Area delle Scienze, 53/a \\
10129 Torino, Italy & Via Trieste 63             & Campus Universitario \\
& 35121 Padova, Italy        & 43100 Parma, Italy \\
\end{tabular}

\bigskip

\leftline{{\it e-mail (DB)}: danilo.bazzanella@polito.it}
\leftline{{\it e-mail (AL)}: languasco@math.unipd.it}
\leftline{{\it e-mail (AZ)}: alessandro.zaccagnini@unipr.it}

\end{document}